\documentstyle[epic,eepic,amssymb]{amsart}
\numberwithin{equation}{section}

\newtheorem {Theorem}            	{Theorem} 
\newtheorem {Lemma}[equation]     	{Lemma}  	
  
\newtheorem {Proposition}[equation]	{Proposition}

\theoremstyle{definition}
\newtheorem{Definition}[equation]{Definition}

\theoremstyle{remark}
\newtheorem*{Remark}{Remark}
\newtheorem{Example}[equation]{Example}
\newtheorem {Corollary} [equation]	{Corollary}  

\def	\DH	{Duistermaat-Heckman\ }
\def	\R	{{\Bbb R}}
\def	\C	{{\Bbb C}}
\def	\CP	{{\Bbb C}{\Bbb P}}
\def	\inv	{^{-1}}
\def	\red	{{\operatorname{red}}}

\def	\tPhi	{\tilde{\Phi}}
\def	\tc	{\tilde{c}}

\def	\g	{{\frak g}}
\def	\t	{{\frak t}}
\def	\h	{{\frak h}}
\def	\to	{{\longrightarrow}}
\def	\ssminus	{{\smallsetminus}}

\newcommand{\labell}[1] {\label{#1}}
\newcommand    {\comment}[1] {}

\begin{document}

\title{Moment maps and non-compact cobordisms}
\author{Yael Karshon}
\email{karshon@@math.huji.ac.il}
\thanks{updated \today}
\thanks{dg-ga/9701006}
\begin{abstract}
We define a moment map for a torus action on a smooth manifold {\em
without\/} a two-form. Cobordisms of such structures are meaningful
even if the manifolds are non compact, as long as the moment maps are
proper. We prove that a compact manifold with a torus action and a
moment map is cobordant the normal bundle of its fixed point set. Two
formulas follow easily: Guillemin's topological version of the abelian
Jeffrey-Kirwan localization, and the Guillemin-Lerman-Sternberg formula
for the \DH measure.
\end{abstract}
\maketitle

\section{Introduction}
This work is part of a joint project with Victor Guillemin and Viktor
Ginzburg, in which we use cobordisms to study group actions on symplectic
manifolds.  Some of our results were announced and their proof sketched
in \cite{GGK}.

In this paper we introduce a new, cleaner, cobordism technique.
We treat noncompact manifolds directly, without 
``approximating" by compact manifolds as we did in \cite{GGK}.
Our basic theorem involves manifolds and not orbifolds; orbifolds only 
come up as reduced spaces.
And we don't need to carry along a symplectic (or presymplectic) form
for the theory to work.


Our basic object is a manifold with a torus action, a proper moment
map, and an equivariant cohomology class.  If $(M,\omega,\Phi)$ is
a symplectic manifold with a moment map in the usual sense, one can
take the equivariant cohomology class represented by the sum $\omega
+ \Phi$.  However, in our new approach
we define a moment map without a two-form (see section
\ref{sec:moment}), and our equivariant cohomology class is unrelated to
the moment map. The normal bundle of every connected component of the
fixed point set is a (non-compact) manifold with the same structure: a
torus action, a proper moment map (that needs to be chosen; see section
\ref{sec:linear}), and an equivariant cohomology class.  We prove in
section \ref{sec:fp} that the original manifold, if compact, is cobordant
to the disjoint union of these normal bundles. The cobording manifold
carries the same kind of structure, and all moment maps are proper.
In particular, if the torus acted with isolated fixed points, the manifold
is cobordant to a disjoint union of vector spaces.

The motivation for this work comes from
the formula of Guillemin-Lerman-Sternberg \cite{GLS}. 
This formula expresses the \DH measure for a compact symplectic
manifold with a torus action with isolated fixed points as the
sum of \DH measures for the tangent spaces at the fixed points.
Viktor Ginzburg observed several years ago \cite{Ginz} that the
\DH measure is an invariant of (compact) cobordism (in the sense of
\cite{Ginz2}). This tempts one to assert that the manifold
is cobordant to the disjoint union of its tangent spaces at the fixed
points. In sections \ref{sec:moment}--\ref{sec:fp} we make sense of this
assertion.  In section \ref{sec:fp} we prove it.  In sections
\ref{sec:reduction} and \ref{sec:cob-red} we state and prove that
``cobordism commutes with reduction", a phenomenon first shown and used
by Guillemin in the context of symplectic manifolds (see \cite{GGK}).
We use this in sections \ref{sec:hamiltonian}, \ref{sec:JK},
and \ref{sec:GLS} to
deduce (generalizations of) the Guillemin-Jeffrey-Kirwan formula and
the Guillemin-Lerman-Sternberg formula.

Another application of our cobordism technique, that will appear
elsewhere (at least in \cite{GGK:book}), is a new proof of the
Atiyah-Bott-Berline-Vergne localization for equivariant differential
forms.  This uses the technique of Lerman cutting applied to spaces with
proper moment maps.

\section*{Acknowledgements}
I thank Victor Guillemin and Viktor Ginzburg for leading to this
work and for useful comments.
Lisa Jeffrey had remarked that in Guillemin's cobordism approach to 
the Jeffrey-Kirwan formula one should be able to replace the moment map
by other maps.

\section{Moment maps}
\labell{sec:moment}
Let $T$ be a torus and $\t$ its Lie algebra.
If $H$ is a closed subgroup of $T$ and $\h$ is the Lie algebra of $H$, 
the {\em $H$-component\/} of a map into $\t^*$ is the composition 
of the map with the natural projection $\t^* \to \h^*$.

\begin{Definition} \labell{def:moment}
Let $M$ a smooth manifold with a smooth $T$-action.
A {\em moment map\/} is a smooth map $ \Phi : M \to \t^* $
with the following properties:
\begin{enumerate}
\item
$\Phi$ is $T$-invariant
\item
For any subgroup $H$ of $T$, the $H$-component of $\Phi$ is constant 
on each connected component of the fixed point set of $H$ in $M$.
\end{enumerate}
\end{Definition}

\begin{Remark}
In property 2 it is enough to take circle subgroups $H$ of $T$. 
\end{Remark}

\begin{Example}
A constant function is a moment map.
\end{Example}

\begin{Remark}
A linear combination of moment maps is a moment map.
Later on we will demand that our moment maps be {\em proper}.
Linear combinations of proper moments maps might
no longer be proper moment maps. However, if $\eta \in \t$ is 
a vector and $\Phi_1$ and $\Phi_2$ are two moment maps whose
$\eta$-components, $\langle \Phi_i , \eta \rangle$, are proper and
bounded from below, then so is every positive combination
of $\Phi_1$ and $\Phi_2$.
\end{Remark}

\begin{Example}
If $T$ is a circle, and we identify $\t^*$ with $\R$,
a moment map is any real valued function that is constant
on each connected component of the set of fixed points.
\end{Example}

\begin{Example}
On a vector space on which the circle
acts linearly and fixes only the origin,
the function $\Phi(v) = ||v||^2$,
with respect to any invariant metric, is a moment map
that is proper and bounded from below.
\end{Example}

\begin{Example}
If $\xi_M$ is a vector field on $M$
that generates a circle action, the function
$\Phi(x) = \langle \xi_x , \xi_x \rangle$, 
where $\langle, \rangle$ is any invariant Riemannian metric,
is a moment map that is proper and bounded from below.
\end{Example}

\begin{Example} \labell{symplectic}
Suppose $(M,\omega)$ is symplectic and $\Phi$ is a moment map
in the traditional sense, i.e.,
$\Phi$ satisfies $\iota(\xi_M) \omega = - d \langle \Phi,\xi \rangle$
for all $\xi \in \t$, where $\xi_M$ is the vector field on $M$
that generates the action of the one parameter subgroup 
$\exp(s\xi)$, $s \in \R$.
Then $\Phi$ is a moment map in the sense of Definition \ref{def:moment}.
\end{Example}

\section{Equivariant cohomlogy}
\labell{sec:eq-coh}
We will work with triples $(M,\Phi,c)$ where $M$ is an oriented manifold
acted upon by a torus $T$, where
$\Phi : M \to \t^*$ is a proper moment map,
and where $c$ is an equivariant cohomology class\footnote{
our cohomology is always over $\R$
} on $M$ (see below).
We keep track of the torus action on $M$ 
although it is not explicit in our notation. 

Recall, an equivariant differential form is a polynomial on 
the Lie algebra $\t$ whose values are differential forms on $M$. 
The equivariant exterior derivative is defined by 
$(d_T \alpha)(\xi) = d(\alpha(\xi)) + \iota(\xi_M) \alpha(\xi)$
for all $\xi \in \t$, $\xi_M$ being the corresponding vector
field on $M$. The equivariant cohomlogy (over $\R$) of a manifold 
is the quotient $\text{ker} (d_T) / \text{image} (d_T)$.
Good references on equivariant cohomology and equivariant
differential forms are \cite{AB} and \cite{BGV}.

For example, if $(M,\omega)$ is a symplectic manifold with a
torus action and $\Phi$ is a moment map, the sum $\omega + \Phi$
is a closed equivariant form. However,
an important point is that {\em one can work with equivariant forms 
that have nothing to do with the moment map}.

\section{Cobordism}
\labell{sec:cobordism}
For cobordism to be meaningful one must make some compactness
assumption. Otherwise, any manifold $M$ would be cobordant to zero
via the noncompact cobordism $M \times (0,1]$.
However, we can drop the compactness assumption if we demand
the cobordism to carry a proper moment map:

\begin{Definition} \labell{def:cobordism}
Let $M, M'$ be smooth oriented manifolds with actions of a torus $T$,
let $\Phi : M \to \t^*$ and $\Phi' : M \to \t^*$ be proper moment maps, 
and let $c$ and $c'$ be equivariant cohomlogy classes on $M$ and $M'$
respectively.
A {\em cobordism\/} between the triples $(M,\Phi,c)$ and $(M',\Phi',c')$ 
is a triple $(W,\tPhi,\tc)$, where $W$ is an oriented 
manifold-with-boundary with a $T$-action, 
$\tPhi : W \to \t^*$ is a proper moment map, 
and $\tc$ is an equivariant cohomology class on $W$, 
and an equivariant diffeomorphism of the boundary $\partial W$ 
with the disjoint union $M \sqcup M'$ that carries $\tPhi$ to 
$\Phi \sqcup \Phi'$, carries $\tc$ to $c \sqcup c'$,
and carries the boundary orientation on $\partial W$
to the given orientation of $M$ and the opposite orientation on $M'$.
\end{Definition}

\begin{Example}
Every ordinary cobordism of compact oriented manifolds
is also a cobordism in the sense of Definition \ref{def:cobordism}
when we take the trivial $T$-actions, zero moment maps,
and zero cohomology classes.
\end{Example}

\begin{Example}
Suppose that $M$ is a manifold with a circle action and that
$\Phi_i : M \to \R$, $i=0,1$, are moment maps
that are proper and bounded from below and that take the same values
on the fixed points.
Then $(M,\Phi_0)$ is cobordant to $(M,\Phi_1)$;
take the cobording manifold $M \times [0,1]$
with the moment map $(m,t) \mapsto (1-t) \Phi_0(m) + t \Phi_1(m)$. 
\end{Example}

\begin{Example}
A space with a trivial torus action and a proper moment map
must be compact, because its moment map must be constant.
Two such spaces are cobordant in the sense of Definition
\ref{def:cobordism} if and only if they are cobordant in the
usual (compact) sense. This is because on a connected 
manifold-with-boundary, if a torus action acts trivially on the boundary
(or, in fact, on any submanifold of codimension one), it acts
trivially everywhere.
\end{Example}

\begin{Remark}
If two compact manifolds with torus actions and proper moment maps
are cobordant in the sense of Definition \ref{def:cobordism} via a
non-compact cobordism, they are also cobordant in the sense of Definition
\ref{def:cobordism} via a compact cobordism.  Thus our notion of cobordism
refines the usual notion of equivariant cobordism for compact manifolds
with torus actions.  To show this, first use Lerman's cutting procedure
\cite{L} to replace a noncompact cobordism by a compact orbifold
cobordism, then get rid of the orbifold singularities by successive
blowups along singular strata (with respect to a complex structure on
the normal bundles to the strata). These blowups are $T$-equivariant,
and the cohomology class pulls back via the blowups.  Details will appear
in \cite{GGK:book}.  In this paper we do not use this result.
\end{Remark}

\begin{Remark}
The above Examples show that our notion
of cobordism is non-trivial. Note that it is absolutely 
crucial for the moment maps to be proper; otherwise,
every manifold $M$ would be cobordant to the empty set
via the non-compact cobordism $M \times (0,1]$,
where the $T$ action and moment maps are induced from 
those on $M$.

Also notice that it is crucial to assume Condition 2 
in Definition \ref{def:moment} of the moment map. 
If we drop this condition and only work with, say, 
manifolds equipped with proper functions to $\R$,
any compact manifold $M$ would be cobordant to the empty set 
via the non-compact cobordism $M \times (0,1]$ and the 
proper function $\tPhi(m,t) = \Phi(m) + 1/t -1$.
\end{Remark}

\begin{Lemma}
Cobordism is an equivalence relation.
\end{Lemma}

\begin{pf}
For smooth manifolds with no extra structure, this is a standard fact
in differential topology. See, e.g., \cite{BJ}.
It is not hard to adapt the standard proof to manifolds equipped with a
$T$-action, a moment map, and an equivariant cohomology class.
We omit the details; we might include them in \cite{GGK:book}.
\end{pf}

\section{Linear actions}
\labell{sec:linear}
The most useful non-compact manifolds are vector spaces.

Recall that a complex representation of the torus $T$ 
is determined up to isomorphism by its {\em weights}, which are 
elements of $\ell^*$, the integral weight lattice in $\t^*$.  
On the underlying real representation, the weights are defined
only up to sign; we call these elements of $\ell^* / \pm 1$
{\em real weights}.  A representation of a torus on a real 
vector space $V$ is determined up to isomorphism by the dimension 
of the subspace, $F$, of fixed vectors, and by the real weights for 
the torus action on the quotient $V / F$.

Fix a real representation of $T$.  
Let $\{ \pm \alpha_i \}$ be its nonzero real weights.
Choose a ``polarizing" vector $\eta \in \t$, a vector such that
$\langle \alpha_i , \eta \rangle \neq 0$ for all $i$.
The corresponding {\em polarized weights\/} are the elements 
$\alpha_i^\eta \in \{ \alpha_i , -\alpha_i \}$ of $\ell^*$
with the sign chosen so that $\langle \alpha_i^\eta , \eta \rangle >0$. 

The following Lemma shows that a moment map for a linear torus
action is determined up to cobordism by its value 
at the origin and its behavior at infinity, provided that 
some component of the moment map is proper.

\begin{Lemma} \labell{lem:linear}
Let $V$ be a real vector space with a linear action of the torus $T$
that fixes only the origin. 
Let $\pm \alpha_1, \ldots, \pm \alpha_m$ be the real weights for the 
action (without repititions).
Let $V = V_1 \oplus \ldots \oplus V_m$ be the decomposition
of $V$ into invariant subspaces such that $T$ acts on $V_i$
by the real weight $\pm \alpha_i$. Any vector $v \in V$
can be expressed uniquely as a sum, $v = \sum v_i$,
where $v_i$ is in $V_i$.
Choose a vector $\eta \in \t$ such that $\langle \alpha_i ,\eta \rangle
\neq 0$ for all $i$. Let $\alpha_1^\eta , \ldots , \alpha_m^\eta$ be the 
corresponding polarized weights.  Choose any invariant metric on $V$. 
Let $a$ be any element of $\t^*$. Then
$$ 
	\Phi(v) = a + \sum_{i=1}^m ||v_i||^2 \alpha_i^\eta
$$
defines a moment map on $V$ whose value at the origin is $a$
and whose $\eta$-component, $\langle \Phi, \eta \rangle$,
is proper and bounded from below. Any other moment map, $\Phi'$, 
with these properties is cobordant to $\Phi$,
i,e., the triple $(V,\Phi',c)$ is cobordant to the triple
$(V,\Phi,c)$ for any equivariant cohomology class $c$ on $V$.
\end{Lemma}

\begin{pf}
Let $\Phi' : V \to \t^*$ be another moment map whose 
value at the origin, $\Phi'(0)$, is $a$, and whose
$\eta$-component, $\langle \Phi,\eta \rangle$, is proper 
and bounded from below.
Let $I$ be the closed interval $[0,1]$.
On the product $I \times V$, take the $T$-action induced from that on $V$.
The function 
$ \tPhi(t,v) := (1-t) \Phi(v) + t \Phi'(v) $
is proper because its $\eta$-component is proper.
It is a moment map because the fixed point set of any subgroup 
of the torus contains the origin, $0$, and because $\Phi(0) = \Phi'(0)$.
The pair $(I \times V , \tPhi)$ provides a cobordism
between $(V,\Phi)$ and $(V,\Phi')$. An equivariant cohomology class
$c$ on $V$ can be simply pulled back to $I \times V$.
\end{pf}

A similar result holds for vector bundles:

\begin{Lemma} \labell{lem:bundle}
Let $E$ be a real vector bundle over a compact connected manifold $F$.
Suppose that the torus $T$ acts linearly on the fibers of $E$,
fixing only the zero section.
Let $\pm \alpha_i \in \t^*$ be the real weights for this action.
Let $E = E_1 \oplus \ldots \oplus E_m$ be the decomposition of $E$
into sub-vector-bundles such that the torus $T$ acts on the fibers 
of $E_i$ by the real weight $\pm \alpha_i$.
Fix a vector $\eta \in \t$ such that $\langle \alpha_i , \eta \rangle
\neq 0$ for all $i$.  Any element $v$ in the vector bundle $E$
can be expressed uniquely as a sum, $v = \sum v_i$,
where $v_i$ is in $E_i$.
Fix an invariant inner product on the fibers of $E$.
Let $a$ be any element of $\t^*$. Then 
$$ \Phi(v) = a + \sum_{i=1}^n ||v_i||^2 \alpha_i^\eta$$
defines a moment map on $E$ whose value on the zero section is $a$
and whose $\eta$-component is proper and bounded from below. 
Any other moment map with these properties is cobordant to $\Phi$.
\end{Lemma}

\begin{pf}
The proof is completely analogous to the proof of Lemma \ref{lem:linear}:
if $\Phi' : E \to \t^*$ is a moment map that takes the value $a$
on the zero section and whose $\eta$-component is proper and bounded
from below, the product $I \times E$ with the moment map 
$\tPhi(t,v) = (1-t)\Phi(v) + t \Phi'(v)$ 
provides the desired cobordism.
\end{pf}

Moreover,

\begin{Lemma} \labell{lem:fancy}
Let $\pi : E \to F$ be a real vector bundle over a compact connected 
manifold.  Suppose that the torus $T$ acts on $E$ by bundle maps
and that $H \subset T$ is a subgroup whose fixed point set
is exactly the zero section. Let $\pm \alpha_i \in \h^*$
be the real weights for the action of the identity component, $H_0$,
of $H$ on the fibers of $E$.
Fix a vector $\eta \in \h$ such that 
$\langle \alpha_i , \eta \rangle \neq 0$ for all $i$. 
Let $\alpha_i^\eta = \pm \alpha_i$ be the polarized weights,
with $\langle \alpha_i^\eta, \eta \rangle >0$.
For $v \in E$, 
let $v = \sum v_i$ be the decomposition such that $H_0$ acts
on $v_i$ with weight $\alpha_i$.  Let
$$ \Phi_1(v) = \sum_{i=1}^n ||v_i||^2 \alpha_i^\eta $$
be the $H$-moment map on the fibers of $E$
as decribed in the previous lemma.
For every embedding $i : \h^* \hookrightarrow \t^*$
whose composition with the restriction $\t^* \to \h^*$ is the identity,
the composition $i \circ \Phi_1 : E \to \t^*$
is a $T$-moment map.

Let $\Phi_2 : F \to \t^*$ be any $T$-moment map
on the base space, $F$, identified with the zero section of $E$. 
(Since $H$ acts trivially and $F$ is connected, 
the image of $\Phi_2$ is contained in a translation of
the annihilator, $\h^0$, of $\h$ in $\t^*$.)
The composition $\Phi_2 \circ \pi : E \to \t^*$ is also a $T$-moment map, 
because the projection $\pi$ is $T$-equivariant. 

The sum $i \circ \Phi_1 + \Phi_2 \circ \pi$ 
is a $T$-moment map on $E$ whose $\eta$-component is proper
and bounded from below.
Any other $T$-moment map on $E$ whose $\eta$-component is proper
and bounded from below
and whose restriction to the zero section is $\Phi_2$
is cobordant to this sum. 
\end{Lemma}

\begin{pf}
Again, if $\Phi' : E \to \t^*$ is such a moment map,
take $I \times E$ with the moment map 
$\tPhi(t,v) = (1-t)\Phi(v) + t \Phi'(v)$
where $\Phi = i \circ \Phi_1 + \Phi_2 \circ \pi$.
\end{pf}

\section{$T$-fixed points}
\labell{sec:fp}
We now state our basic theorem.

\begin{Theorem} \labell{basic}
Let $M$ be an oriented manifold with an action of a torus $T$, let $\Phi$
be a moment map, and let $c$ be an equivariant cohomology class on $M$.
Suppose that the vector $\eta \in \t$ is such that the $\eta$-component
of the moment map,
$$
	\langle \Phi,\eta \rangle : M \to \R,
$$
is proper and bounded from below.
Denote by $M^\eta$ the zero set of the vector field on $M$
corresponding to $\eta$.
Every connected component, $F$, of $M^\eta$
is a compact $T$-invariant submanifold of $M$.
Let $NF$ be the normal bundle of $F$ in $M$
with the $T$-action induced from that on $M$.
Let $c_F$ be the pullback of the cohomology class $c$ via the maps
$NF \to F \hookrightarrow M$.
Let $\Phi_F^\eta$ be a moment map on $NF$ defined by
\begin{equation} \labell{Phip}
	\Phi_F^\eta(v) = \Phi(p) + \sum_i ||v_i||^2 \alpha_{i,F}^\eta 
 \quad v \in N_pF
\end{equation} 
as in the of notation Lemma \ref{lem:fancy} with $E=NF$.
Then the triple $(M,\Phi,c)$ is cobordant to the disjoint union,
over all the components $F$ of the set $M^\eta$,
of the triples $(NF,\Phi_F^\eta,c_F)$.
\end{Theorem}

An important special case is that of a torus action with 
isolated fixed points on a compact manifold:

\begin{Theorem}
Let $M$ be a compact oriented manifold with an action of a torus $T$
with isolated fixed points, let $\Phi$ be a moment map, and let $c$ be 
an equivariant cohomology class on $M$. 
Let $\eta \in \t$ be any nonzero vector
that does not belong to the Lie algebra of the stablizer 
of any point in $M$.  For every fixed point, $p$, in $M$, 
let $T$ act on the tangent space $T_pM$ by the linear isotropy $T$-action,
and let $c_p$ be the pullback of the equivariant
cohomology class $c$ via the maps $T_pM \to \{p\} \hookrightarrow M$.
Let $\Phi_p^\eta$ be a moment map on $NF$ defined as in Lemma
\ref{lem:linear} with $V = T_pM$. Then the triple $(M,\Phi,c)$
is cobordant to the disjoint union
of the triples $(T_pM ,\Phi_p^\eta,c_p)$
over all fixed points $p$ in $M$.
\end{Theorem}

\begin{pf*}{Proof of a special case of Theorem \ref{basic}}
To convey the idea of the proof let us first deal with 
a circle action with isolated fixed points. Let
$ \Phi : M \to \R$
be a moment map that is proper and bounded from below.
We would like to show that the pair $(M,\Phi)$ is cobordant to the 
union of the pairs $(T_pM,\Phi_p)$ 
over the fixed points, $p$, of the circle action,
where $\Phi_p : T_pM \to \R$ is defined by $\Phi_p(v) = \Phi(p) + ||v||^2$
for some invariant metric on $T_pM$.

Take the product of $M$ with the half-open interval:
$$
	M \times (0,1] \/ ,
$$
with the circle action induced from that on $M$.
This provides a non-compact cobordism of $M$ with the 
empty set. The function $(m,x) \mapsto \Phi(m)$
is a moment map whose restriction to the boundary, 
$M \times \{ 1  \}$, is $\Phi$. However, this moment map is not proper.
We fix this by adding a function $\rho(x)$, $0 < x \leq 1$,  
that approaches $\infty$ as $x$ approaches $0$, 
and that vanishes for $\epsilon \leq x \leq 1$.
If there are no fixed points, the function
$$\tPhi(m,x) = \Phi(m) + \rho(x)$$
is a proper moment map whose restriction to the boundary,
$M \times \{ 1 \}$, is $\Phi$, and we are done.

If there are fixed points, this does not work:
the components of the fixed point set in $M \times (0,1]$ 
have the form $F \times (0,1]$ where $F$ is a component
of the fixed point set in $M$.
The function $\tPhi$ is constant on $F \times [\epsilon,1]$
but not on $F \times (0,\epsilon]$, so it is not constant
on the fixed components, so it is not a moment map. 
However, in this case we can take the
same function $\tPhi$ on the manifold 
$$(M \times (0,1]) \ssminus \sqcup_p B_p \/ ,$$
where $B_p$ is an $\epsilon$-neighborhood of $(p,0)$
in $M \times (0,1]$ with respect to some invariant metric on $M$.
The union is over the (isolated) fixed points $p$ in $M$.
Here $\epsilon$ is the same as in the definition of $\tPhi$,
and the metric on $M$ is chosen so that $\epsilon$-neighborhoods
of fixed points are disjoint from each other, and so that
an $\epsilon$-neighborhood of the fixed point $p$
is equivariantly diffeomorphic to $T_pM$.
The removal of $B_p$ creates a boundary component
that is equivariantly diffeomorphic to a ball in $M$ around $p$
via the projection map $(m,x) \mapsto m$.
This is, further, equivariantly diffeomorphic to $T_pM$.
The function $\tPhi$ transforms to a function $\Phi_p(v)$ on $T_pM$
that takes the value $\Phi(p)$ at the origin 
and that approaches $\infty$ as we approach infinity in $T_pM$.
By Lemma \ref{lem:linear}, this function is cobordant to $\Phi_p$.
\end{pf*}

\begin{pf*}{Proof of Theorem \ref{basic}.}
We return to the full generality of torus actions.
Let $H \subseteq T$ be the closure of the one parameter subgroup
$\{ \exp(s\eta) , s \in \R \}$ of $T$.
Since the set $M^\eta$ coincides with the fixed point set, $M^H$,
of $H$ in $T$, it consists of a disjoint union of submanifolds of $M$.

Let $F$ be a connected component of the set $M^\eta$.
Since the $\eta$-component of $\Phi$ is proper and is constant on $F$,
\ $F$ is compact. 
There exists an invariant Riemannian metric on $M$
and an $\epsilon >0$
such that the $\epsilon$-neighborhood of $F$ in $M$
is equivariantly diffeomorphic to the normal bundle, $NF$, 
of $F$ in $M$, and such that these neighorhoods are disjoint
for different $F$'s. 

Let $B_F$ denote the $\epsilon$-neighborhood
of $\{ 0 \} \times F$ in $[0,1] \times M$.
As before, $W := (M \times (0,1]) \ssminus \sqcup_F B_F$
provides a noncompact equivariant cobordism between 
$M$ and $\sqcup_F NF$. 
Our cobording moment map will now have the form
$$
	\tPhi(m,x) = \Phi(m) + \rho(m,x).
$$

We first define $\rho$ locally:
let $U$ be an invariant neighborhood of $m$ in $M$ 
that admits an equivariant retraction to the orbit $T \cdot m$.
On the product $U \times (0,1]$, define
$$
	\rho_m(m',x) = \alpha_m w(x) \/ ,
$$
where $w : (0,1] \to \R$ is such that $w(x)$
approaches infinity as $x$ approaches $0$ and vanishes
for $ \epsilon \leq x \leq 1$,
and where $\alpha_m \in \t^*$ is chosen in the following way.
If $m \in M^\eta$ we set $\alpha_m =0$.
If $m \not \in M^\eta$ then the vector $\eta$ is not in the Lie algebra,
$\g_m$, of the stabilizer of $m$.  We then choose $\alpha_m \in \t^*$
to be an element of the annihilator, $\g_m^0$, of $\g_m$, such that 
$\langle \alpha,\eta \rangle >0$. 

We patch the functions $\rho_m$ together, using a partition of unity,
to obtain a moment map $\rho$ on 
$W = M \times (0,1] \ssminus \sqcup_F B_F$
satisfying $\langle \rho(m,x),\eta \rangle \to \infty$ as $x \to 0$
and $\rho(m,x) = 0$ for $\epsilon \leq x \leq 1$.

Finally, take
$$
	\tPhi(m,x) := \Phi(m) + \rho(m,x) .
$$
This is a moment map on $(M \times (0,1]) \ssminus \sqcup_F B_F$.
The component $\langle \tPhi,\eta \rangle$ is proper, 
so $\tPhi$ itself is proper.
On the boundary component $M \times \{ 1 \}$ we get the original moment
map:
$$
	\tPhi(m,1) = \Phi(m).
$$
The new boundary component that we get by removing the open set $B_F$ 
can be identified with an $\epsilon$-neighborhood of $F$
via $(m,x) \mapsto m$, and further identified with the normal
bundle $NF$. Via these identifications, the moment map on this
boundary component becomes a moment map $\Phi_F'$ on $NF$ 
that coincides with the map $\Phi$ at the points of the zero section
(identified with the points of the base $F$),
and whose $\eta$-component, $\langle \Phi_F'(v),\eta \rangle$, 
approaches infinity as $v$ approaches infinity.
By Lemma \ref{lem:fancy}, this is further cobordant to 
the moment map $\Phi_F^\eta$ 
that is given by the formula \eqref{Phip}.

Finally, the pullback of the equivariant cohomology class $c$
under the projection $M \times (0,1] \to M$
restricts to an equivariant cohomology class $\tc$ on the
cobording manifold $W$, whose pullback to the boundary 
component $\partial B_F \cong N_F$ is $c_F$.
\end{pf*}

\section{Reduction}
\labell{sec:reduction}
If $(M,\omega)$ is a symplectic manifold with a torus action
and $\Phi$ is a moment map, the torus acts locally freely 
on the regular level sets of $\Phi$.
The same holds for moment maps in the sense of Definition \ref{def:moment};
the proof is only slightly harder:

\begin{Lemma} \labell{orbifold}
Let $M$ be a manifold with an action of a torus $T$, and let $\Phi$ be a 
moment map.  If $a$ is a regular value of $\Phi$, the $T$-action
on the level set $\Phi\inv(a)$ is locally free,
i.e., all stabilizers are discrete.
\end{Lemma}

\begin{pf}
First let us assume that $T$ is a circle.  Let $F$ be a component 
of its fixed point set. 
Since the moment map is constant on $F$, the restriction of $d\Phi$
to the tangent bundle of $F$ is zero.
Since the moment map is invariant under the circle action,
the restriction of $d\Phi$ to the normal bundle of $F$ is zero.
This shows that every fixed point for the circle action is a critical
point for the moment map.

Now let the torus $T$ be of any dimension.
We need show that any point $m \in M$ whose stabilizer has positive 
dimension is a critical point for $\Phi$.
Let $H$ be a circle subgroup of the stabilizer of $m$.
The previous paragraph, applied to the action of $H$,
implies that $m$ is a critical point for the $H$-component
of $\Phi$, hence for $\Phi$.
\end{pf}

We denote by $M_\red = M_\red^\Phi(a)$ the quotient $\Phi\inv(a) / T$.
By Lemma \ref{orbifold}, if $a$ is a regular value of $\Phi$, 
$M_\red$ is an orbifold. If $\Phi$ is proper, $M_\red$ is compact.

An orientation on $M$ induces an orientation on $M_\red$ in the 
following way.
At any point $p \in \Phi\inv(a)$ we have a (non-canonical) splitting:
$T_p M = T_{\Phi(p)} M_\red \oplus \t \oplus \t^*$ 
where the $\t$-piece is the tangent to the orbit and where
the $\t^*$-piece is the normal to the level set $\Phi\inv(a)$, 
identified with $\t^*$ via $d\Phi$.
We demand that the orientation of $M_\red$, followed by the
natural orientation on $\t \oplus \t^*$, give the orientation of $M$. 

Consider the inclusion map $i$ and the quotient map $\pi$
in the following diagram:
\begin{equation} \labell{diagram}
\begin{array}{ccc}
 \Phi\inv(a) & \stackrel{i}{\hookrightarrow} & M \\
 \pi \downarrow \phantom{\pi} & & \\
 M_\red & &
\end{array}
\end{equation}
An equivariant cohomlogy class on $M$ descends to an ordinary 
cohomology class on $M_\red$ in the following way.
An ordinary cohomology class on $M_\red$ can be viewed
as an equivariant cohomology class
when $M_\red$ is taken with the trivial circle action.
As such, the class $c_\red$ on $M_\red$ that is associated
to an equivariant class $c$ on $M$ is characterized by the equality
\begin{equation} \labell{cartan}
	\pi^* (c_\red) = i^* (c)
\end{equation}
of equivariant classes on the level set $\Phi\inv(a)$.
Here we used the fact that if a group acts locally
freely on a space, the equivariant cohomology (over $\R$) 
is equal to the ordinary cohomology of the quotient.

The map $c \mapsto c_\red$ that we have just described is a ring
homomorphism from the equivariant cohomology of $M$ to the ordinary
cohomology of $M_\red$.  This is the famous {\em Kirwan map}.

\begin{Definition}
Let $M$ be a manifold with an action of a torus $T$, let $\Phi$
be a moment map on $M$, and let $a$ be a regular value of $\Phi$.
The {\em reduction\/} of the triple $(M,\Phi,c)$ at the value $a$ 
is the pair $(M_\red,c_\red)$
where $M_\red$ is the compact oriented orbifold $\Phi\inv(a) / T$
and where $c_\red$ is the ordinary (real) cohomology class  
on $M_\red$ satisfying \eqref{cartan}.
\end{Definition}

To avoid ambiguity we may keep track of $\Phi$ and $a$
in the notation; we may write $M_\red(a)$ and $c_\red(a)$,
or $M_\red^\Phi(a)$ and $c_\red^\Phi(a)$.

\section{Cobordism of reduced spaces}
\labell{sec:cob-red}
Reduction of a triple $(M,\Phi,c)$, consisting of
a manifold with $T$-action, proper moment map, and 
equivariant cohomology class, produced a pair
$(M_\red,c_\red)$ where $M_\red$ is a compact orbifold and $c_\red$
is an ordinary cohomology class on $M_\red$.

\begin{Definition}
A cobordism between pairs $(M,c)$ and $(M',c')$,
where $M$ and $M'$ are oriented orbifolds and $c$ and $c'$ are cohomology
classes, is a pair $(W,\tc)$ where $W$ is an oriented
orbifold-with-boundary and $\tc$ an ordinary cohomology class on $W$,
and a diffeomorphism of $M \sqcup M'$ with the boundary 
$\partial W$, that sends $\tc$ to $c \sqcup c'$,
that respects the orientation on $M$, and that flips the orientation
on $M'$.
\end{Definition}

The integral $ \int_M c$ is an invariant of cobordism of the pair 
$(M,c)$; this follows from Stokes's theorem.

\begin{Lemma} \labell{regular}
Let $(M,\Phi,c)$ be a manifold with a $T$-action,
a proper moment map, and an equivariant cohomology class.
Let $a$ and $b$ be two values in $\t^*$
such that all values in the closed interval $[a,b]$
are regular values for $\Phi$. 
Then $M_\red(a)$ is cobordant to $M_\red(b)$.
\end{Lemma}

\begin{pf}
The quotient $\Phi\inv([a,b]) / T$ provides the desired cobordism.
The cobording cohomology class is the one whose pullback 
to $\Phi\inv([a,b])$ is the restriction of $c$.
\end{pf}

\begin{Remark}
Under the assumptions of Lemma \ref{regular}, the reduced spaces 
are diffeomorphic and not just cobordant. 
(Without the torus action and cohomology class, this follows from
Ehresmann's fibration theorem.)
\end{Remark}

Lemma \ref{regular} shows that, up to cobordism, the reduced space 
depends in a robust way on the value at which we reduce. 
In section \ref{sec:fp} we will see that the dependence 
on the moment map is also robust.

Guillemin was the first to note that ``cobordism commutes with 
reduction" is true and is useful for symplectic manifolds. 
Here is the statement for moment maps in our new sense:

\begin{Lemma} \labell{cob-red}
Let $(M,\Phi,c)$ and $(M',\Phi',c')$ be cobordant oriented manifolds 
with $T$-actions, proper moment maps, and equivariant cohomology classes.
Let $a \in \t^*$ be a value that is regular for both $\Phi$ and $\Phi'$.
Then the corresponding reductions, $(M_\red,c_\red)$ and 
$(M'_\red , c'_\red)$, are also cobordant.
\end{Lemma}

\begin{pf}
Let $(W,\tPhi,\tc)$ be a cobording manifold, proper moment map,
and equivariant cohomology class.
Since $\Phi$ and $\Phi'$ are proper, their sets of regular values 
are open, so any value $b$ that is close enough to $a$ 
is also a regular value for $\Phi$ and for $\Phi'$.
For such $b$, the reduced spaces at $b$ are
cobordant to the reduced spaces at $a$, by Lemma \ref{regular}.
If we choose such a $b$ that is also a regular value for $\tPhi$, 
the quotient $\tPhi\inv(b) / T$ provides a cobordism
between the reduced spaces at $b$. 
The cobording cohomology class is the one whose pullback 
to $\tPhi\inv(b)$ is the restriction of $\tc$.
\end{pf}

\begin{Corollary} \labell{invt}
For every regular value $a \in \t^*$, the integral 
$\int_{M_\red(a)} c_\red(a)$ 
is an invariant of cobordism of the triple $(M,\Phi,c)$.
\end{Corollary}

\section{Hamiltonian manifolds}
\labell{sec:hamiltonian}
Corollary \ref{invt} implies that the \DH measure of
a Hamiltonian $T$-manifold is an invariant of cobordism.
We first recall some definitions:

\begin{Definition}
Let $(M,\omega)$ and $(M',\omega')$ be two symplectic manifolds 
acted upon by a torus $T$ with proper moment maps $\Phi$ and $\Phi'$.
A {\em cobordism\/} between them
is an oriented manifold-with-boundary $W$ with a $T$-action,
a closed two-form $\tilde{\omega}$, and a proper moment map
$\tilde{\Phi} : W \to \t^*$ ({\it i.e.}, a map satisfying 
$d \langle \tilde{\Phi}, \xi \rangle = \iota(\xi_W) \tilde{\omega}$
for all $\xi \in \t$,
when $\xi_W$ is the vector field on $W$ corresponding to $\xi$),
and an equivariant diffeomorphism of $M \sqcup M'$ with the boundary
$\partial W$ that carries $\tilde{\omega}$ to $\omega \sqcup \omega'$
and $\tilde{\Phi}$ to $\Phi \sqcup \Phi'$,
and that respects orientation on $M$ and flips orientation on $M'$.
\end{Definition}

\begin{Remark}
In the above definition there is no reason other
than psychological to assume that the $\omega$ and $\omega'$
are nondegenerate. They can be just closed two-forms. 
In fact, any two closed two-forms on the same manifold that belong to the 
same cohomology class are cobordant, so we can even work with
cohomology classes intead of with two-forms.
\end{Remark}

\begin{Definition}
{\em Liouville measure\/} on a symplectic manifold $(M,\omega)$ is
defined on open sets
by integrating the volume form $\omega^n / n!$ \ ($n = \dim M$)
with respect to the symplectic orientation.  
In the presence of a $T$-action, the corresponding {\em \DH measure\/} 
is the measure on $\t^*$ obtained as the push-forward of
Liouville measure by the moment map.  
\end{Definition}

Let the torus $T$ act on a symplectic manifold $(M,\omega)$ 
with a proper moment map $\Phi$. 
At any regular value $a$, the Marsden-Weistein reduced space
is the compact symplectic orbifold $(M_\red,\omega_\red)$
defined by $M_\red = \Phi \inv(a) / T$
and $\pi^* \omega_\red = i^* \omega$,
where $i$ and $\pi$ are the inclusion and quotient maps
as in \eqref{diagram}.

\begin{Definition}
The {\em \DH function\/} associates to a regular value 
of the moment map the symplectic volume of the reduced space: 
$a \mapsto \int_{M_\red} \omega_\red^d / d!$, \ 
$M_\red = \Phi\inv(a)/T$, \  $d = \dim M_\red$.  
\end{Definition}

The \DH measure is absolutely continuous with respect to
Lebesgue measure on $\t^*$, provided that the torus action is effective.
It is well known, and is not hard to prove, that the Radon-Nikodym 
derivative of the \DH measure is the \DH function.

The Kirwan map (see \S\ref{sec:reduction}) sends the equivariant 
cohomology class $[\omega + \Phi]$ on $M$ to the ordinary cohomology 
class $[\omega_\red + a]$ (of mixed degree) on $M_\red(a)$.

Finally, Corollary \ref{invt} applied to the equivariant cohomology 
classes of $(\omega + \Phi - a)^d / d!$ on $M$ implies the following:

\begin{Proposition} \labell{DH}
The \DH measure is an invariant of cobordism. 
\end{Proposition}

\begin{Remark}
The fact that the \DH measure is invariant
under (compact) cobordisms was observed by V. Ginzburg several
years ago. This was one of the main motivations for introducing 
cobordism techniques into equivariant symplectic geometry.
\end{Remark}

\section{Application: the Jeffrey-Kirwan localization formula}
\labell{sec:JK}

Theorem \ref{basic}, together with the facts that 
``cobordism commutes with reduction" (Lemma \ref{cob-red}) and that 
the integral of a reduced class is an invariant of cobordism
(Corollary \ref{invt}), yields a formula for integrals over
the reduced space in terms of fixed point data:

\begin{Proposition} \labell{JKG}
Let $M$ be an oriented manifold with an action of a torus, $T$.
Let $\Phi$ be a moment map, and let $c$ be an equivariant
cohomology class on $M$. Suppose that the vector $\eta \in \t$
is such that the $\eta$-component of the moment map,
$$ \langle \Phi,\eta \rangle : M \to \R,$$
is proper and bounded from below. Denote by $M^\eta$
the zero set of the vector field on $M$ corresonding to $\eta$.
For every connected component, $F$, of $M^\eta$, denote by 
$NF$ the normal bundle of $F$ in $M$, and let
$\Phi_F^\eta : NF \to \t^*$
the proper moment described in Formula \eqref{Phip}.
Let $c_F$ be the equivariant cohomology class on $NF$ obtained
by pulling back $c$ via the inclusion $F \hookrightarrow M$
and further via the bundle map $NF \to F$.
Let $a \in \t^*$ be a value that is regular for $\Phi$ and for
all the $\Phi^\eta_F$'s. 
Let $M_\red = M_\red^\Phi(a)$ and 
$(NF)_\red = (NF)_\red^{\Phi_F^\eta}(a)$
be the corresponding reduced spaces, and
let $c_\red = c_\red^\Phi(a)$ and 
$(c_F)_\red = (c_F)_\red^{\Phi^\eta_F}(a)$
be the corresponding cohomology classes on these reduced spaces.
Then
\begin{equation} \labell{JKG-formula}
\int_{M_\red} c_\red = 
  \sum_F \int_{NF_\red} (c_F)_\red.
\end{equation}
\end{Proposition}

Notice that in \eqref{JKG-formula}, 
if we choose $\eta$ generic, the sum is over the components 
of the fixed point set, $M^T$. Also notice that the $F$-summand
is nonzero if and only if $\langle \Phi(F), \eta \rangle < 0$.

Proposition \ref{JKG} generalizes Guillemin's topological form
for the Jeffrey-Kirwan localization theorem in the abelian case.
Jeffrey and Kirwan consider Hamiltonian actions
of compact nonabelian groups on compact symplectic manifolds.
Every cohomology class on the reduced space at $0$ 
comes via ``Kirwan's map" from an equivariant cohomology class
on the manifold (as in \eqref{cartan}).  
Jeffrey and Kirwan express the integral of such a class over the 
reduced space by an explicit formula that only depends on data 
at the set of points that are fixed under the action of a maximal torus.
Shaun Martin showed, by topological means, that their formula
would follow from a localization formula for abelian groups. 
Guillemin interpreted this abelian localization formula
as Formula \eqref{JKG-formula}.

\section{Application: the Guillemin-Lerman-Sternberg formula}
\labell{sec:GLS}

A special case of Proposition \ref{JKG} 
is the formula of Guillemin, Lerman, and Sternberg
for the \DH measure. (The definition of the \DH 
measure was recalled in \S\ref{sec:hamiltonian}.)

\begin{Proposition}
Let $(M,\omega)$ be a presymplectic oriented manifold 
with a Hamiltonian action of the torus $T$ 
and with a moment map $\Phi : M \to \t^*$
(in the traditional sense:  $d \langle \Phi, \xi \rangle
 = - \iota (\xi_M) \omega$ for all $\xi \in \t$, where $\xi_M$
is the corresponding vector field on $M$).
Let $\eta \in \t$ be a vector such that the function
$\langle \Phi , \eta \rangle : M \to \R$ is proper and bounded from below.
Let $M^\eta$ be the set of zeros of the vector field on $M$
that corresponds to $\eta$.

For every component $F$ of the set $M^\eta$,
let $NF$ be the normal bundle of $F$ in $M$,
choose an invariant metric on the fibers of $NF$,
and let $\Phi_F^\eta$ be the moment map described in Theorem \ref{basic}.
There exists a symplectic form on the fibers of $NF$
for which $\Phi_F^\eta$, restricted to each fiber, is a moment map
in the usual, symplectic, sense.
With these forms $NF$ becomes a symplectic vector bundle. 
Let $\omega_F$ be any closed two-form on the total space 
of the bundle $NF$ whose restriction to the fibers is given 
by the symplectic forms just mentioned, and whose
pullback to $F$ via the zero section
coincides with the pullback of $\omega$ under the inclusion
of $F$ into $M$.
(Such a closed two-form, $\omega_F$, can be constructed by ``minimal 
coupling"; see, e.g., \cite{GLS:book}.  Notice that $\omega_F$ need not 
be symplectic on $NF$.) 
Then $\Phi_F^\eta : NF \to \t^*$ is a moment map for $(NF,\omega_F)$.

Denote by $\text{DH}(M)$ the \DH measure corresponding to the triple
$(M,\omega,\Phi)$, and by $\text{DH}(NF)$ the \DH measure corresponding 
to the triple $(NF,\omega_F,\Phi_F)$. Then 
\begin{equation} \labell{GLS-formula}
 \text{DH}(M) = \sum_F \text{DH}(NF) . 
\end{equation}
\end{Proposition}

The Guillemin-Lerman-Sternberg formula 
was originally asserted for compact symplectic manifolds
and for torus actions with isolated fixed points.
The vector $\eta \in \t$ was chosen to be generic,
in the sense that the zero set $M^\eta$ was required to be equal to $M^T$, 
the fixed point set for $T$ in $M$. 
The right hand side of \eqref{GLS-formula} in this case is the sum 
of the \DH measures, $\text{DH}(T_pM)$, for the linear torus actions
on the tangent spaces at the fixed point.
The original formula of Guillemin-Lerman-Sternberg \cite{GLS}
involved an explicit formula for the measure $\text{DH}(T_pM)$
as the push-forward of Lebesgue measure on a positive orthant in $\R^n$, 
\ $n = {1 \over 2} \dim M$, via a linear projection to $\t^*$. 
See \cite{GLS:book}.

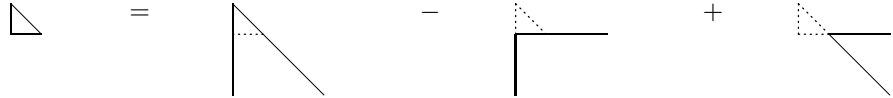
\begin{figure}
\setlength{\unitlength}{0.008in}
\hfill
\begin{picture}(20,60)
\put(0,60) {\line(0,-1){20}}
\put(0,60) {\line(1,-1){20}}
\put(0,40) {\line(1,0){20}}
\end{picture}
\hfill
	\begin{picture}(30,60)(0,0)
	\put(10,50){$=$}
	\end{picture}
\hfill
\begin{picture}(60,60)
\put(0,60) {\line(0,-1){60}}
\put(0,60) {\line(1,-1){60}}
\dottedline{4}(0,40)(20,40)
\end{picture}
\hfill
	\begin{picture}(30,60)(0,0)
	\put(15,50){$-$}
	\end{picture}
\hfill
\begin{picture}(60,60)(0,0)
\dottedline{4}(0,60)(0,40)
\dottedline{4}(0,60)(20,40)
\put(0,40) {\line(1,0){60}}
\put(0,40) {\line(0,-1){40}}
\end{picture}
\hfill
	\begin{picture}(30,60)(0,0)
	\put(15,50){$+$}
	\end{picture}
\hfill
\begin{picture}(60,60)(0,0)
\dottedline{4}(0,60)(0,40)
\dottedline{4}(0,60)(20,40)
\dottedline{4}(0,40)(20,40)
\put(20,40) {\line(1,0){40}}
\put(20,40) {\line(1,-1){40}}
\end{picture}
\hfill
\caption{The Guillemin-Lerman-Sternberg formula for $\CP^2$}
\label{triangle}
\end{figure}

\begin{Example}
\labell{ex:CP2}
Take $\CP^2$ with the Fubini-Study symplectic form,
whose pullback to $S^5 \subset \C^3$ 
is the standard two-form $\sum_{j=1}^3 dx_j \wedge dy_j$.
The linear action of $T = S^1 \times S^1$ on $\C^3$
by scalar multiplication on the first two coordinates
induces an action on $\CP^2$.
The image of the moment map is Lebesgue measure on a triangle. 
The action on $\CP^2$ has three isolated fixed points; their images
are the vertices of the triangle. 
Each of the summands $\text{DH}(T_pM)$ is 
Lebesgue measure on the region between two rays and zero outside.
The Guillemin-Lerman-Sternberg formula exhibits the triangle as a 
combination of three such ``wedges"; see figure \ref{triangle}.
Notice that this involves a choice of direction
in which the infinite rays are pointing; this is the choice
of the vector $\eta$. 
\end{Example}

\begin{figure}
\setlength{\unitlength}{0.01in}
\hfill
\begin{picture}(20,60)
\put(0,60) {\line(0,-1){20}}
\put(0,60) {\line(1,-1){20}}
\put(0,40) {\line(1,0){20}}
\end{picture}
\hfill
	\begin{picture}(30,60)(0,0)
	\put(10,50){$=$}
	\end{picture}
\hfill
\begin{picture}(60,60)
\put(0,60) {\line(0,-1){20}}
\put(0,60) {\line(1,-1){60}}
\put(0,40) {\line(1,0){60}}
\put(2,42) {$+$}
\put(50,20) {$-$}
\end{picture}
\hfill
	\begin{picture}(30,60)(0,0)
	\put(15,50){$+$}
	\end{picture}
\hfill
\begin{picture}(60,60)(0,0)
\dottedline{4}(0,60)(0,40)
\dottedline{4}(0,60)(20,40)
\dottedline{4}(0,40)(20,40)
\put(20,40) {\line(1,0){40}}
\put(20,40) {\line(1,-1){40}}
\end{picture}
\hfill
\caption{The (generalized) Guillemin-Lerman-Sternberg formula for $\CP^2$,
with a non-generic choice of ``polarizing vector"}
\label{triangle2}
\end{figure}
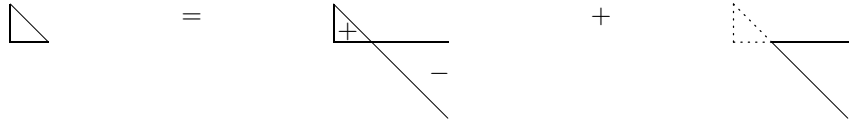

\begin{Example}
In Example \ref{ex:CP2}, if $\eta$ is chosen non-generically, 
the zero set $M^\eta$
consists of one copy of $\CP^1$ and one isolated fixed point.
The (generalized) Guillemin-Lerman-Sternberg formula then gives 
the combination illustrated in Figure \ref{triangle2}.
Here the first summand is the \DH measure for the normal bundle
of $\CP^1$ in $\CP^2$. This summand is a {\em signed\/} measure; 
this reflect the fact that the cohomology class on this normal bundle 
does not have an invariant symplectic representative.
\end{Example}

\begin{thebibliography}{GGK2}

\bibitem[AB]{AB} M. F. Atiyah and R. Bott,
  {\em The moment map and equivariant cohomology},
  Topology {\bf 23} (1984), 1--28.

\bibitem[BGV]{BGV} N. Berline, E. Getzler, and M. Vergne,
  {\em Heat Kernels and Dirac Operators},
  Springer Verlag, 1992.

\bibitem[BJ]{BJ} Th.\ Br\"ocker and K. J\"anich,
  {\em Introduction to Differential Topology},
  Cambridge University Press, 1982.

\bibitem[G1]{Ginz} V. L. Ginzburg, private communication, December 1993.

\bibitem[G2]{Ginz2} V. L. Ginzburg, 
{\em Calculation of contact and symplectic cobordism groups},
Topology {\bf 31} (1992), 767--773.

\bibitem[GGK1]{GGK}
  V. L. Ginzburg, V. Guillemin, and Y. Karshon,
  {\em Cobordism theory and localization formulas for Hamiltonian
  group actions}, Int.\ Math.\ Res.\ Notices {\bf 5} (1996), 221--234.

\bibitem[GGK2]{GGK:book}
  V. L. Ginzburg, V. Guillemin, and Y. Karshon,
  monograph to be published by the AMS University Lecture Series,
  in preparation.

\bibitem[GLS1]{GLS}
  V. Guillemin, E. Lerman, and S. Sternberg, {\em On the Kostant
  multiplicity formula}, J.\ Geom.\ Phys.\ {\bf 5} (1988), no.\ 4,
  721--750.

\bibitem[GLS2]{GLS:book}
  V. Guillemin, E. Lerman, and S. Sternberg, {\em Symplectic
  Fibrations and Multiplicity Diagrams}, Cambridge Univ.\ Press,
  1996.

\bibitem[L]{L}
  E. Lerman, {\em Symplectic Cuts},
  Math.\ Res.\ Letters {\bf 2} (1995), 247--258.

\end{thebibliography}
\end{document}